\documentclass[11pt]{amsart}
\def\stackunder#1#2{\mathrel{\mathop{#2}\limits_{#1}}}
\newcommand{\C}{{\if mm {{\rm C}\mkern -15mu{\phantom{\rm t}\vrule}}
\mkern +10mu \else \leavemode \hbox{I}\kern -.17em \hbox{C} \fi}}
\newcommand{\Q}{{\if mm {{\rm Q}\mkern -16mu{\phantom{\rm t}\vrule}}
\mkern +10mu \else \leavemode \hbox{I}\kern -.17em \hbox{Q} \fi}}
\newcommand{\R}{{\if mm {\rm I}\mkern -3mu{\rm R}\else \leavevmode
\hbox{I}\kern -.17em\hbox{R} \fi}}
\newcommand{\N}{{\if mm {\rm I}\mkern -3mu{\rm N}\else \leavevmode
\hbox{I}\kern -.17em \hbox{N} \fi}}

\begin{document}
\title{Homotopical structures in categories}
\author{Adrian Duma}
\address{Department of Mathematics,
University of Craiova,
13 A.I. Cuza, 1100 Craiova}
\email{ady@royal.net}
\author{Cristian Vladmirescu}
\email{vladimirescucris@hotmail.com}
\date{}

\begin{abstract}
In this paper is presented a new approach to the axiomatic homotopy theory
in categories, which offers a simpler and more useful answer to this old
question: how two objects in a category (without any topological feature)
can be deformed each in other?
\end{abstract}
\keywords{Category, Homotopical structure}
\subjclass[1991]{46M20, 46M99, 46A55}
\maketitle

\section{Introduction}

As we know, homotopy theory is a branch of algebraic topology concerned with
the study of those properties of topological spaces that are invariant under
homotopy equivalence. Homotopy theoretic methods can be applied to a variety
of interesting problems, for example to the classification of manifolds
(simply connected $5-$manifolds are understood and classified by homotopy
theory) and in establishing important fixed point theorems.

In the forty-five years since Eilenberg and Steenrod wrote down the axioms
for homology theory there have been repeated attempts to ''do the same
thing'' homotopy theory. Various categorical settings for homotopical
algebra have been proposed, among which we briefly recall the following 
ones.

a) Kan [8] developed an ''abstract homotopy theory'' in the framework of
categories equipped with a cylinder functor; other extensions of these
approaches are in Kamps [7], Hubner [6] and Kleisli [9].

b) Keller [5] introduced the definition of a homotopy theory as a
hyperfunctor, taking small categories, functors and natural transformations
to categories, functors and natural transformations and enjoying a list of
properties which thus serve as the axioms for homotopy theory.

It is the purpose of this paper to introduce a new approach to the axiomatic
homotopy theory in categories, which offers an useful answer to the question
if two objects in a category (without any topological feature) could be
deformed each in other.

Our viewpoint is that the categorical foundations of homotopy theory can be
based on very simple categorical notions. Namely, a homotopy theory is a way
of associating with each object of a category another object and two
morphisms in the same category, satisfying four easy axioms. Brutally
summarized, these axioms say that the homotopy between morphisms is a
congruence relation (i.e. an equivalence relation which is also functorial).

In the second section we give the definition of a homotopical structure
which generates an internal homotopy theory in an arbitrary category which,
in turn, gives rise to a congruence relation between morphisms.

In Section 3 there are given some examples and applications. Firstly, the
classical homotopy theory of the category of topological spaces. Next, there
are constructed homotopy theories for convex compact sets in locally convex
spaces and for Banach spaces.

\section{Homotopical structures}

Let ${\C}$ be a category. In what follows, $\left| {\C}\right| $
will denote the class of the objects in ${\C}$ and, if $X,$ $Y\in \left|
{\C}\right| $, then we shall write ${\C}\left( X,Y\right) $ for the
set of all morphisms from $X$ to $Y.$

{\bf Definition 1. }{\it A {\bf homotopical structure} on the category }$%
{\C}$ {\it is a correspondence rule }
\[
S\stackrel{H}{\rightarrow }\left( \hat{S};i_{S},j_{S}\right)
\]
{\it which associates to each object }$S${\it \ from }${\C}$ {\it %
another object }$\hat{S}${\it \ and two morphisms }$i_{S},$ $j_{S}\in {\C%
}\left( S,\hat{S}\right) $ {\it with the following properties}:

I. {\it For every }$S\in \left| {\C}\right| ,$ {\it there exists }$p\in
{\C}\left( \hat{S},S\right) $, {\it such that }$pi_{S}=pj_{S}=1_{S}$;

II. {\it For every }$S\in \left| {\C}\right| ,$ {\it there exists }$k\in
{\C}\left( \hat{S},\hat{S}\right) $, {\it such that }$ki_{S}=j_{S}$ {\it %
and }$kj_{S}=i_{S}$;

III. {\it For every }$S,$ $T\in \left| {\C}\right| $ {\it and every }$h,$
$h^{\ast }\in {\C}\left( \hat{S},T\right) $ {\it such that }$h^{\ast
}i_{S}=hj_{S},$ {\it there exists }$h^{\ast \ast }\in {\C}\left( \hat{S}%
,T\right) $, {\it such that }$hi_{S}=h^{\ast \ast }i_{S}$ {\it and }$h^{\ast
}j_{S}=h^{\ast \ast }j_{S}$;

IV. {\it For every }$S,$ $T\in \left| {\C}\right| $ {\it and every }$%
u\in {\C}\left( S,T\right) ,$ {\it there exists }$\hat{u}\in {\C}%
\left( \hat{S},\hat{T}\right) $, {\it such that }$\hat{u}i_{S}=i_{T}u$ {\it 
and} $\hat{u}j_{S}=j_{T}u$.

{\bf Definition 2. }{\it If }$R,$ $Q\in \left| {\C}\right| $ {\it and }$%
\varphi ,$ $\psi \in {\C}\left( R,Q\right) $ {\it we say that the
morphisms }$\varphi ${\it \ and }$\psi ${\it \ are }${\C-}${\bf homotopic%
}{\it \ iff there exists }$h\in {\C}\left( \hat{R},Q\right) ${\it \ such
that }$hi_{R}=\varphi $ {\it and }$hj_{R}=\psi .$

We shall write $\varphi \stackrel{H}{\sim }\psi .$

{\bf Theorem 1. }{\it The binary relation ''}$\stackrel{H}{\sim }${\it ''
attached to a homotopical structure is a congruence relation.}

{\bf Proof. }We denote by $\stackunder{X,Y}{\stackrel{H}{\sim }}$ the binary
relation induced on the set ${\C}\left( X,Y\right) $.

A. The relation $\stackunder{X,Y}{\stackrel{H}{\sim }}$ is an equivalence
relation.

A$_{1}.$ The reflexivity: $\left( \forall \right) $ $u\in {\C}\left(
X,Y\right) $ we have $u\stackunder{X,Y}{\stackrel{H}{\sim }}u.$ We claim
that $\left( \exists \right) $ $h\in {\C}\left( \hat{X},Y\right) $ such
that $hi_{X}=u$ and $hj_{X}=u.$ It suffices to choose $h=up$, where $p\in
{\C}\left( \hat{X},X\right) $ is given by the axiom I.

A$_{2}.$ The symmetry: $\left( \forall \right) $ $u,$ $v\in {\C}\left(
X,Y\right) $ such that $u\stackunder{X,Y}{\stackrel{H}{\sim }}v$ we also
have $v\stackunder{X,Y}{\stackrel{H}{\sim }}u.$ There exists $h\in {\C}%
\left( \hat{X},Y\right) $ with $hi_{X}=u,$ $hj_{X}=v.$ Let $h^{\ast }=hk$, 
$%
k\in {\C}\left( \hat{X},\hat{X}\right) $ given by axiom II. We have
\[
h^{\ast }i_{X}=hki_{X}=hj_{X}=v
\]
and
\[
h^{\ast }j_{X}=hi_{X}=u,
\]
thus $v\stackunder{X,Y}{\stackrel{H}{\sim }}u.$

A$_{3}.$ The transitivity: $\left( \forall \right) $ $u,$ $v,$ $w\in {\C}%
\left( X,Y\right) $ such that $u\stackunder{X,Y}{\stackrel{H}{\sim }}v$ $\ 
$%
and $v\stackunder{X,Y}{\stackrel{H}{\sim }}w$ we also have 
$u\stackunder{X,Y%
}{\stackrel{H}{\sim }}w.$ There exist $h$ and $h^{\ast }\in C\left( \hat{X}%
,Y\right) $ such that
\begin{eqnarray}
hi_{X} &=&u, \\
hj_{X} &=&v, \\
h^{\ast }i_{X} &=&v, \\
h^{\ast }j_{X} &=&w.
\end{eqnarray}

From $\left( 2\right) ,$ $\left( 3\right) $ and axiom III it follows that
there is $h^{\ast \ast }\in C\left( \hat{X},Y\right) $ such that
\begin{equation}
h^{\ast \ast }i_{X}=hi_{X}
\end{equation}
\qquad\ and
\begin{equation}
h^{\ast \ast }j_{X}=hj_{X}.
\end{equation}
\qquad

From $\left( 5\right) $ and $\left( 1\right) $ we obtain $h^{\ast \ast
}i_{X}=u$ and from $\left( 6\right) $ and $\left( 4\right) $ we obtain $%
h^{\ast \ast }j_{X}=w.$ Hence, $u\stackunder{X,Y}{\stackrel{H}{\sim }}w.$

B. The compatibility of the binary relation ''$\stackrel{H}{\sim }$'' with
the composition of morphisms.

Let $u\in {\C}\left( X,Y\right) $ be such that $v\stackunder{X,Y}{%
\stackrel{H}{\sim }}v.$ Let $X^{\prime },$ $Y^{\prime }\in \left| {\C}%
\right| $ and the morphisms $f\in {\C}\left( X^{\prime },X\right) $ and $%
g\in {\C}\left( Y,Y^{\prime }\right) .$ We claim that
\[
guf\stackunder{X^{\prime },Y^{\prime }}{\stackrel{H}{\sim }}gvf.
\]

There exists $h\in {\C}\left( \hat{X},Y\right) $ such that $hi_{X}=u,$ $%
hj_{X}=v.$ We have
\[
ghi_{X}f=guf
\]
and
\[
ghj_{X}f=gvf.
\]

By considering the objects $X^{\prime },$ $X\in \left| {\C}\right| $ and
the morphism $f\in {\C}\left( X^{\prime },X\right) $ and using axiom IV
we obtain the existence of a morphism $\hat{f}\in {\C}\left( \hat{X}%
^{\prime },\hat{X}\right) $ such that $\hat{f}i_{X^{\prime }}=i_{X}f$ and $%
\hat{f}j_{X^{\prime }}=j_{X}f.$ Then we define $h^{\ast }\in {\C}\left(
\hat{X}^{\prime },Y^{\prime }\right) $ by $h^{\ast }=gh\hat{f}$. Then we
have
\[
h^{\ast }i_{X^{\prime }}=gh\hat{f}i_{X^{\prime }}=ghi_{X}f=guf
\]
and
\[
h^{\ast }j_{X^{\prime }}=gh\hat{f}j_{X^{\prime }}=ghj_{X}f=gvf,
\]
thus
\[
guf\stackunder{X^{\prime },Y^{\prime }}{\stackrel{H}{\sim }}gvf.
\]

\hfill $\Box $

\section{Examples and applications}

{\bf 3.1.} 
${\C}=$ Top, $H$ is the classical homotopical relation. When $%
S\in \left| \mbox{Top}\right| $ we have
\begin{eqnarray*}
\hat{S} &=&S\times \left[ 0,1\right] ; \\
i_{S},\mbox{ }j_{S} &\in &\mbox{ Top }\left( S,S\times \left[ 0,1\right]
\right) , \\
i_{S}\left( s\right)  &=&\left( s;0\right) , \\
j_{S}\left( s\right)  &=&\left( s;1\right) ,\mbox{ for every }s\in S.
\end{eqnarray*}

{\bf 3.2.} Let ${\Q}$ be the category where the objects are compact
convex subsets of Hausdorff locally convex spaces and, if $S,$ $T\in \left|
{\Q}\right| $, then the set of morphisms from $S$ to $T$ is:
\[
{\Q}\left( S,T\right) =\left\{ f:S\rightarrow T,\mbox{ }f\mbox{ affine,
continuous, }f\left( \mbox{ex }S\right) \subseteq \mbox{ex }T\right\} .
\]

Here ex $K$ denotes the set of all extremal points of the compact convex 
$K.$

We introduce the following notations:
\[
P=\left\{ \mu \in C\left( \left[ 0,1\right] \right) ^{\ast },\mbox{ }\mu
\geq 0,\mbox{ }\left\| \mu \right\| =1\right\} ;
\]
$\Delta =$ ex $P$, the family of Dirac measures; $S\hat{\otimes}_{\pi }P,$
the projective tensor product of $S$ and $P$, where $S\in \left| {\Q}%
\right| $. In what follows,
\[
\omega _{\pi ,S}:S\times P\rightarrow S\hat{\otimes}_{\pi }P
\]
denotes the canonical biaffine continuous map. We define
\[
i_{S},j_{S}\in {\Q}\left( S,S\hat{\otimes}_{\pi }P\right)
\]
by
\[
\left.
\begin{array}{c}
i_{S}\left( s\right) =\omega _{\pi ,S}\left( s;\delta _{0}\right) \\
j_{S}\left( s\right) =\omega _{\pi ,S}\left( s;\delta _{1}\right)
\end{array}
\right. ,\mbox{ for every }s\in S,
\]
where $\delta _{0},$ $\delta _{1}\in \Delta $ are the Dirac measures
attached to $0$ and $1,$ respectively.

Let us remark that $i_{S}$ and $j_{S}$ are affine, continuous and, clearly, 
$%
s\in $ ex $S$ implies that
\[
\left\{ \left( i_{S}\left( s\right) ,j_{S}\left( s\right) \right) \right\}
\subset \omega _{\pi ,S}\left( \mbox{ex }\left( S\times P\right) \right) =%
\mbox{ ex }\left( S\hat{\otimes}_{\pi }P\right) .
\]

{\bf Theorem 2.}{\it \ The correspondence}
\[
S\rightarrow \left( S\hat{\otimes}_{\pi }P;i_{S},j_{S}\right)
\]
{\it defines a homotopical structure on the category }${\Q}$.

{\bf Proof. }

I. Let $s\in \left| {\Q}\right| $. We introduce the continuous biaffine
mapping $\hat{p}:S\times P\rightarrow S,$
\[
\hat{p}\left( s,\mu \right) =s,\mbox{ for every }s\in S,\mbox{ }\mu \in P.
\]

Then there exists a unique affine continuous mapping $p:S\hat{\otimes}_{\pi
}P\rightarrow S$ such that
\[
p\omega _{\pi ,S}=\hat{p}.
\]

We have
\[
p\left( \mbox{ex }\left( S\hat{\otimes}_{\pi }P\right) \right) =p\omega
_{\pi ,S}\left( \mbox{ex }\left( S\times P\right) \right) =\hat{p}\left(
\mbox{ex }\left( S\times P\right) \right) =\mbox{ ex }S
\]
and thus $p\in {\Q}\left( S\hat{\otimes}_{\pi }P,S\right) .$ Let us
remark that, if $s\in S$, then we have
\[
pi_{S}\left( s\right) =p\omega _{\pi ,S}\left( s;\delta _{0}\right) 
=\hat{p}%
\left( s;\delta _{0}\right) =s
\]
and
\[
pj_{S}\left( s\right) =p\omega _{\pi ,S}\left( s;\delta _{1}\right) 
=\hat{p}%
\left( s;\delta _{1}\right) =s.
\]

It follows that $pi_{S}=pj_{S}=1_{S}.$

II. Let us consider the isometry $w\in L\left( C\left( \left[ 0,1\right]
\right) ,C\left( \left[ 0,1\right] \right) \right) ,$
\[
w\left( u\right) \left( t\right) =u\left( 1-t\right) ,\mbox{ for all }u\in
C\left( \left[ 0,1\right] \right) \mbox{ and }t\in \left[ 0,1\right] .
\]

Obviously, $w^{\ast }\left( P\right) =P$ and $w^{\ast }\left( \delta
_{\lambda }\right) =\delta _{1-\lambda }$, for every $\lambda \in \left[ 
0,1%
\right] .$ For every $S\in \left| {\Q}\right| $ we define the biaffine
continuous mapping $k_{0}:S\times P\rightarrow S\hat{\otimes}_{\pi }P$,
expressed by
\[
k_{0}\left( s,\mu \right) =\omega _{\pi ,S}\left( s;w^{\ast }\left( \mu
\right) \right) .
\]

Then, there exists a unique affine continuous mapping
\[
k:S\hat{\otimes}_{\pi }P\rightarrow S\hat{\otimes}_{\pi }P
\]
such that
\[
k\omega _{\pi ,S}=k_{0}.
\]

If $x\in $ ex $\left( S\hat{\otimes}_{\pi }P\right) $, then there exist $%
s\in $ ex $S$ and $\lambda \in \left[ 0,1\right] $ such that
\[
x=\omega _{\pi ,S}\left( s;\delta _{\lambda }\right) .
\]

Then we have
\begin{eqnarray*}
k\left( x\right) &=&k\omega _{\pi ,S}\left( s;\delta _{\lambda }\right)
=k_{0}\left( s;\delta _{\lambda }\right) =\omega _{\pi ,S}\left( s;w^{\ast
}\left( \delta _{\lambda }\right) \right) = \\
&=&\omega _{\pi ,S}\left( s;\delta _{1-\lambda }\right) \in \omega _{\pi
,S}\left( \mbox{ex }\left( S\times P\right) \right) =\mbox{ ex }\left( 
S\hat{%
\otimes}_{\pi }P\right) .
\end{eqnarray*}

It follows that $k\in {\Q}\left( S\hat{\otimes}_{\pi }P,S\hat{\otimes}%
_{\pi }P\right) $. On the other part, if $s\in S$, then we obtain
\begin{eqnarray*}
ki_{S}\left( s\right)  &=&k\omega _{\pi ,S}\left( s;\delta _{0}\right)
=k_{0}\left( s;\delta _{0}\right) =\omega _{\pi ,S}\left( s;w^{\ast }\left(
\delta _{0}\right) \right) = \\
&=&\omega _{\pi ,S}\left( s;\delta _{1}\right) =j_{S}\left( s\right)
\end{eqnarray*}
and, analogously, $kj_{S}=i_{S}.$

III. Let $S,$ $T\in \left| {\Q}\right| $ and $h,$ $h^{\ast }\in {\Q}%
\left( S\hat{\otimes}_{\pi }P,T\right) $ be such that
\[
h^{\ast }i_{S}=hj_{S}.
\]

We obtain that
\[
h\omega _{\pi ,S}\left( s;\delta _{1}\right) =h^{\ast }\omega _{\pi
,S}\left( s;\delta _{0}\right) ,\mbox{ for every }s\in S.
\]

We define the mapping $h_{0}:S\times \Delta \rightarrow T$, expressed by
\begin{eqnarray*}
h_{0}\left( s;\delta _{\lambda }\right) &=&h\omega _{\pi ,S}\left( s;\delta
_{2\lambda }\right) ,\mbox{ for every }s\in S,\mbox{ }\lambda \in \left[ 0,%
\frac{1}{2}\right] \\
h_{0}\left( s;\delta _{\lambda }\right) &=&h^{\ast }\omega _{\pi ,S}\left(
s;\delta _{2\lambda -1}\right) ,\mbox{ for every }s\in S,\mbox{ }\lambda \in 
\left[ \frac{1}{2},1\right] .
\end{eqnarray*}

Evidently, $h_{0}$ is continuous and the map $h_{0}\left( \cdot ,\delta
_{\lambda }\right) $ is affine for every $\lambda \in \left[ 0,1\right] $.
Using Krein-Milman Theorem, it follows that $h_{0}$ admits a unique
continuous biaffine extension $h_{1}:S\times P\rightarrow T$. Then, there
exists a unique continuous affine mapping
\[
h^{\ast \ast }:S\hat{\otimes}_{\pi }P\rightarrow T,
\]
such that $h^{\ast \ast }\omega _{\pi ,S}=h_{1}.$

Let $x\in $ ex $\left( S\hat{\otimes}_{\pi }P\right) $ and $s\in $ ex $S$, 
$%
\gamma \in \left[ 0,1\right] $ be such that
\[
x=\omega _{\pi ,S}\left( s;\delta _{\gamma }\right) .
\]

If $\gamma \in \left[ 0,\frac{1}{2}\right] $, we have
\begin{eqnarray*}
h^{\ast \ast }\left( x\right)  &=&h^{\ast \ast }\omega _{\pi ,S}\left(
s;\delta _{\gamma }\right) =h_{1}\left( s;\delta _{\gamma }\right)
=h_{0}\left( s;\delta _{\gamma }\right) = \\
&=&h\omega _{\pi ,S}\left( s;\delta _{2\gamma }\right) \in h\left( \omega
_{\pi ,S}\left( \mbox{ex }\left( S\times P\right) \right) \right) = \\
&=&h\left( \mbox{ex }\left( S\hat{\otimes}_{\pi }P\right) \right) \subseteq
\mbox{ ex }T
\end{eqnarray*}
and the case $\gamma \in \left[ \frac{1}{2},1\right] $ is completely
similar. Thus, we obtain that
\[
h^{\ast \ast }\in {\Q}\left( S\hat{\otimes}_{\pi }P,T\right) .
\]

When $s\in S$, we have
\begin{eqnarray*}
h^{\ast \ast }i_{S}\left( s\right)  &=&h^{\ast \ast }\omega _{\pi ,S}\left(
s;\delta _{0}\right) =h_{1}\left( s;\delta _{0}\right) =h_{0}\left( s;\delta
_{0}\right) = \\
&=&h\omega _{\pi ,S}\left( s;\delta _{0}\right) =hi_{S}\left( s\right) ,
\end{eqnarray*}
hence
\[
h^{\ast \ast }i_{S}=hi_{S}
\]
and, analogously,
\[
h^{\ast \ast }j_{S}=h^{\ast }j_{S}.
\]

IV. Let $S,$ $T\in \left| {\Q}\right| $ and let $u\in {\Q}\left(
S,T\right) $. We define the biaffine continuous map
\[
u_{0}:S\times P\rightarrow T\hat{\otimes}_{\pi }P
\]
by
\[
u_{0}\left( s;\mu \right) =\omega _{\pi ,T}\left( u\left( s\right) ,\mu
\right) ,
\]
where $s\in S$ and $\mu \in P.$ Then there is a (unique) continuous affine
map $\hat{u}:S\hat{\otimes}_{\pi }P\rightarrow T\hat{\otimes}_{\pi }P$, such
that $\hat{u}\omega _{\pi ,S}=u_{0}.$ When $x\in $ ex $\left( 
S\hat{\otimes}%
_{\pi }P\right) $, there exist $s\in $ ex $S$ and $\lambda \in \left[ 0,1%
\right] $ such that
\[
x=\omega _{\pi ,S}\left( s;\delta _{\gamma }\right) .
\]

Let us remark that
\[
u\left( s\right) \in \mbox{ ex }T.
\]

We have
\begin{eqnarray*}
\hat{u}\left( x\right)  &=&\hat{u}\omega _{\pi ,S}\left( s;\delta _{\gamma
}\right) =u_{0}\left( s;\delta _{\gamma }\right) =\omega _{\pi ,T}\left(
u\left( s\right) ;\delta _{\gamma }\right) \in  \\
&\in &\omega _{\pi ,T}\left( \mbox{ex }\left( T\times P\right) \right) =%
\mbox{ ex }\left( T\hat{\otimes}_{\pi }P\right) .
\end{eqnarray*}

Thus, $\hat{u}\in {\Q}\left( S\hat{\otimes}_{\pi }P,T\hat{\otimes}_{\pi
}P\right) .$ We obtain
\[
\hat{u}i_{S}\left( s\right) =\hat{u}\omega _{\pi ,S}\left( s;\delta
_{0}\right) =u_{0}\left( s;\delta _{0}\right) =\omega _{\pi ,T}\left(
u\left( s\right) ;\delta _{0}\right) =i_{T}u\left( s\right) ,
\]
thus
\[
\hat{u}i_{S}=i_{T}u
\]
and, clearly,
\[
\hat{u}j_{S}=j_{T}u.
\]

\hfill $\Box $

{\bf 3.3. Homotopical Banach spaces}

Let $SB$ be the category in which the objects are Banach spaces and the
morphisms are linear continuous operators.

For each $E\in \left| SB\right| $, we denote by $\varphi _{E}\in L\left(
E,E^{\prime \prime }\right) $ the canonical map. For $E\in \left| SB\right| 
$
we define $\hat{E}\in \left| SB\right| $ and the morphisms $i_{E},$ $%
j_{E}\in L\left( E,\hat{E}\right) $ by
\[
\hat{E}=\left( E\times E^{\prime \prime }\times E^{\prime \prime }\right)
_{2},\mbox{ }\left.
\begin{array}{c}
i_{E}\left( x\right) =\left( x,\varphi _{E}\left( x\right) ,0\right)  \\
j_{E}\left( x\right) =\left( x,0,\varphi _{E}\left( x\right) \right)
\end{array}
\right. ,\mbox{ for every }x\in E.
\]

{\bf Theorem 3.}{\it \ The correspondence}
\[
E\stackrel{H}{\rightarrow }\left( \hat{E};i_{E},j_{E}\right)
\]
{\it defines a homotopical structure on the category }$SB.$

{\bf Proof.}

I. Let $E\in \left| SB\right| $. We define $p\in L\left( \hat{E},E\right) $
by $p\left( x,x^{\prime \prime },y^{\prime \prime }\right) =x.$ Then,
clearly, $pi_{E}\left( x\right) =pj_{E}\left( x\right) =x.$

II. If $E\in \left| SB\right| $, we introduce the morphism $k\in L\left(
\hat{E},\hat{E}\right) $ given by
\[
k\left( x,x^{\prime \prime },y^{\prime \prime }\right) =\left( x,y^{\prime
\prime },x^{\prime \prime }\right) .
\]

Then we have $ki_{E}=j_{E}$ and $kj_{E}=i_{E}.$

III. Let $E,$ $F\in \left| SB\right| $ and $h,$ $h^{\ast \ast }\in L\left(
\hat{E},F\right) $ such that $h^{\ast }i_{E}=hj_{E}.$ We define $h^{\ast
\ast }\in L\left( \hat{E},F\right) $, expressed by
\[
h^{\ast \ast }\left( x,x^{\prime \prime },y^{\prime \prime }\right) =h\left(
x,x^{\prime \prime },y^{\prime \prime }\right) +h^{\ast }\left( x,x^{\prime
\prime },y^{\prime \prime }\right) -h\left( x,0,\varphi _{E}\left( x\right)
\right) .
\]

Then we have
\begin{eqnarray*}
h^{\ast \ast }i_{E}\left( x\right)  &=&h\left( x,\varphi _{E}\left( x\right)
,0\right) +h^{\ast }\left( x,\varphi _{E}\left( x\right) ,0\right) -h\left(
x,0,\varphi _{E}\left( x\right) \right) = \\
&=&hi_{E}\left( x\right) +h\left( x,0,\varphi _{E}\left( x\right) \right)
-h\left( x,0,\varphi _{E}\left( x\right) \right) =hi_{E}\left( x\right) , \\
h^{\ast \ast }j_{E}\left( x\right)  &=&h\left( x,0,\varphi _{E}\left(
x\right) \right) +h^{\ast }\left( x,0,\varphi _{E}\left( x\right) \right)
-h\left( x,0,\varphi _{E}\left( x\right) \right) = \\
&=&h^{\ast }j_{E}\left( x\right) .
\end{eqnarray*}

IV. Let $E,$ $F\in \left| SB\right| $ and $u\in L\left( E,F\right) $. We
define $\hat{u}\in L\left( \hat{E},\hat{F}\right) ,$ by
\[
\hat{u}\left( x,x^{\prime \prime },y^{\prime \prime }\right) =\left(
ux,u^{\prime \prime }x^{\prime \prime },u^{\prime \prime }y^{\prime \prime
}\right) .
\]

Then we have
\[
\hat{u}i_{E}\left( x\right) =\hat{u}\left( x,\varphi _{E}\left( x\right)
,0\right) =\left( ux,u^{\prime \prime }\varphi _{E}\left( x\right) ,0\right)
=\left( ux,\varphi _{F}\left( ux\right) ,0\right) =i_{F}u\left( x\right) .
\]

Analogously we obtain $\hat{u}j_{E}=j_{F}u.$\hfill $\Box $

{\bf Theorem 4. }{\it Let }$E,$ $F\in \left| SB\right| $ {\it and }$u,$ $%
v\in L\left( E,F\right) $. {\it Then }$u\stackrel{H}{\sim }v$ {\it if and
only if there exists }$w\in L\left( E^{\prime \prime },F\right) $ {\it such
that }$w\varphi _{E}=u-v.$

{\bf Proof. }Suppose that $u\stackrel{H}{\sim }v$ and let $h\in L\left( 
\hat{%
E},F\right) $ be such that $hi_{E}=u$ and $hj_{E}=v.$ We define $w\in
L\left( E^{\prime \prime },F\right) $ by
\[
w\left( x^{\prime \prime }\right) =h\left( 0,x^{\prime \prime },-x^{\prime
\prime }\right) ;
\]
then
\begin{eqnarray*}
w\varphi _{E}\left( x\right) &=&h\left( 0,\varphi _{E}\left( x\right)
,-\varphi _{E}\left( x\right) \right) =h\left( x,\varphi _{E}\left( x\right)
,0\right) -h\left( x,0,\varphi _{E}\left( x\right) \right) = \\
&=&u\left( x\right) -v\left( x\right) .
\end{eqnarray*}

Conversely, if $w\varphi _{E}=u-v$, then we define $h\in L\left( \hat{E}%
,F\right) $ by
\[
h\left( x,x^{\prime \prime },y^{\prime \prime }\right) =vx+wx^{\prime \prime
}.
\]

Then we obtain
\[
hi_{E}\left( x\right) =h\left( x,\varphi _{E}\left( x\right) ,0\right)
=vx+w\varphi _{E}\left( x\right) =vx+ux-vx=ux
\]
and
\[
hj_{E}\left( x\right) =h\left( x,0,\varphi _{E}\left( x\right) \right) =vx.
\]

{\bf Definition 3.} {\it The Banach spaces }$E${\it \ and }$F${\it \ are
called {\bf homotopic} iff there exists the operators }$u\in L\left(
E,F\right) ${\it \ and }$v\in L\left( F,E\right) ${\it \ such that}
\[
vu\stackrel{H}{\sim }1_{E}\mbox{{\it \mbox{ }and\mbox{ }}}uv\stackrel{H}{%
\sim }1_{F}.
\]

{\bf Definition 4.} {\it The Banach space }$E${\it \ is called {\bf %
contractible} iff the spaces }$E$ {\it and }$\left\{ 0\right\} ${\it \ are
homotopic}.

{\bf Theorem 5.}{\it \ The Banach space }$E$ {\it is contractible if and
only if its canonical image in }$E^{\prime \prime }$ {\it is complemented.}

{\bf Proof.} Suppose that the Banach space $E$ is contractible. Then, using
Theorem 4, there is $w\in L\left( E^{\prime \prime },E\right) $ such that $%
w\varphi _{E}=1_{E}.$ Now, set $p=\varphi _{E}w.$ It follows that $p$ is a
projector from all the space $E^{\prime \prime }$ onto the canonical image
of $E.$

Conversely, let $p:E^{\prime \prime }\rightarrow \varphi _{E}\left( E\right)
$ be a projector. We introduce the operator
\[
\varphi _{E}^{-1}:\varphi _{E}\left( E\right) \rightarrow E
\]
and define the operator $w$ by
\[
w=\varphi _{E}^{-1}p.
\]

Then, $w\in L\left( E^{\prime \prime },E\right) $ and we have $w\varphi
_{E}=1_{E}.$ Using Theorem 4, it follows that $1_{E}\stackrel{H}{\sim }0$
and thus $E$ is contractible.\hfill $\Box $

Therefore, we have the following classification:

\begin{quotation}
{\it Contractible spaces: }reflexive Banach spaces, $KB-$spaces, $l_{\infty
},$ $L^{\infty }\left( \left[ 0,1\right] \right) ;$

{\it Non-contractible spaces: }the preduals of $l_{1}.$
\end{quotation}

{\bf Example.}

The morphisms $1_{c_{0}}$ and $0$ are not homotopic. Indeed, suppose that
there is $w\in L\left( l_{\infty },c_{0}\right) $ such that $w\varphi
_{c_{0}}=1_{c_{0}}.$ Then, $c_{0}$, being separable, it follows that $w$ is
weakly compact. Consequently, $1_{c_{0}}$ must be weakly compact, which, in
turn, implies the reflexivity of $c_{0},$ contradiction.

{\bf Proposition 1.}{\it \ Let }$E,$ $F\in \left| SB\right| $ {\it and let 
}$%
u,$ $v\in L\left( E,F\right) $ {\it be such that the operator }$u-v$ {\it is
weakly compact. Then }$u\stackrel{H}{\sim }v.$

{\bf Proof.} Let $w=u-v.$ Since $w$ is weakly compact, we have $w^{\prime
\prime }\left( E^{\prime \prime }\right) \subset F$ and, consequently, $%
w=w^{\prime \prime }\varphi _{E}.$ Using Theorem 4, we obtain that $u%
\stackrel{H}{\sim }v.$\hfill $\Box $

{\bf Proposition 2.}{\it \ Let }$E$ {\it be a Banach space with separable
bidual and let }$F=c_{0}${\it . Then any two morphisms }$u,$ $v\in L\left(
E,F\right) $ {\it are homotopic}$.$

{\bf Proof.} Let $w=u-v.$ Since $c_{0}$ has the separable extension
property, there exists $\widetilde{w}\in L\left( E^{\prime \prime },F\right)
$, such that $\widetilde{w}\varphi _{E}=w.$ Therefore, we obtain that $u%
\stackrel{H}{\sim }v,$ using Theorem 4.\hfill $\Box $

{\bf Proposition 3. }{\it \ Let }$E$ $=c_{0}$ {\it and let }$F$ {\it be a
separable Banach space. If }$u,$ $v\in L\left( E,F\right) $, {\it then }$u$
{\it and }$v$ {\it are homotopic if and only if the operator }$u-v$ {\it is
weakly compact.}

{\bf Proof. }If $u$ and $v$ are homotopic, then, by setting $w=u-v$, there
is $f\in L\left( l_{\infty },F\right) $ such that $f\varphi _{c_{0}}=w.$
Since $F$ is separable, it follows that $f$ is weakly compact. Therefore, 
$w$
is weakly compact, too. The converse follows from Proposition 1.\hfill $\Box$

\end{document}